\documentclass[12pt]{article}
\usepackage{amssymb}
\usepackage{amsmath}
\usepackage{amsthm}
\usepackage{color}
\usepackage{comment}
\usepackage{geometry}		
\usepackage{graphicx}

\let\originalleft\left
\let\originalright\right
\renewcommand{\left}{\mathopen{}\mathclose\bgroup\originalleft}
\renewcommand{\right}{\aftergroup\egroup\originalright}

\begin{document}

\def\co{o}
\def\cO{O}
\def\ee{\varepsilon}
\def\manualEndProof{\hfill$\Box$}
\def\rD{{\rm D}}
\def\re{{\rm e}}
\def\ri{{\rm i}}

\newcommand{\myStep}[2]{{\bf Step #1} --- #2}

\newcommand{\removableFootnote}[1]{}

\newtheorem{theorem}{Theorem}[section]
\newtheorem{corollary}[theorem]{Corollary}
\newtheorem{lemma}[theorem]{Lemma}
\newtheorem{proposition}[theorem]{Proposition}

\theoremstyle{definition}
\newtheorem{definition}{Definition}[section]
\newtheorem{example}[definition]{Example}

\theoremstyle{remark}
\newtheorem{remark}{Remark}[section]



\title{
Hopf-like boundary equilibrium bifurcations involving two foci in Filippov systems.
}
\author{
D.J.W.~Simpson\\\\
Institute of Fundamental Sciences\\
Massey University\\
Palmerston North\\
New Zealand
}
\maketitle


\begin{abstract}

This paper concerns two-dimensional Filippov systems --- ordinary differential equations
that are discontinuous on one-dimensional switching manifolds.
In the situation that a stable focus transitions to an unstable focus
by colliding with a switching manifold as parameters are varied,
a simple sufficient condition for a unique local limit cycle to be created is established.
If this condition is violated, three nested limit cycles may be created simultaneously.
The result is achieved by constructing a Poincar\'e map
and generalising analytical arguments that have been employed for continuous systems.
Necessary and sufficient conditions for the existence of pseudo-equilibria
(equilibria of sliding motion on the switching manifold) are also determined.
For simplicity only piecewise-linear systems are considered.

\end{abstract}

\section{Introduction}
\label{sec:intro}
\setcounter{equation}{0}

Physical systems involving impacts, switches, thresholds and other abrupt events
are often well modelled by ordinary differential equations that are piecewise-smooth.
The phase space of a piecewise-smooth system contains switching manifolds 
where the functional form of the equations changes.
As parameters are varied an equilibrium may collide with a switching manifold ---
this is known as a boundary equilibrium bifurcation (BEB) \cite{DiBu08}.
There are many possibilities for the dynamics near a BEB,
including chaos in systems of three or more dimensions \cite{Gl18,Si16c,Si18d}.

This paper concerns two-dimensional systems of the form
\begin{equation}
\begin{bmatrix} \dot{x} \\ \dot{y} \end{bmatrix} =
\begin{cases} F_L(x,y;\mu), & x < 0, \\ F_R(x,y;\mu), & x > 0, \end{cases}
\label{eq:ode}
\end{equation}
where $F_L$ and $F_R$ are smooth vector fields, and $\mu \in \mathbb{R}$ is a parameter.
The system \eqref{eq:ode} has the single switching manifold $x=0$.

If \eqref{eq:ode} is continuous on $x=0$ (i.e.~$F_L(0,y;\mu) \equiv F_R(0,y;\mu)$),
then, assuming genericity conditions are satisfied,
BEBs involve two equilibria (one for each of $F_L$ and $F_R$).
These coincide at the BEB and a limit cycle is created in some cases \cite{FrPo98,Si10}.
If instead \eqref{eq:ode} is discontinuous on $x=0$ (i.e.~\eqref{eq:ode} is a Filippov system \cite{DiBu08,Fi88}),
then orbits may slide on $x=0$.
Generic BEBs involve one equilibrium and one pseudo-equilibrium (an equilibrium of the sliding vector field).
As in the continuous case, these coincide at the BEB and a limit cycle may be created \cite{Gl16d,HoHo16,KuRi03}.

BEBs can mimic classical bifurcations, such as saddle-node bifurcations and Hopf bifurcations \cite{LeNi04}.
In particular, if \eqref{eq:ode} is continuous on $x=0$
and the BEB involves an unstable focus for $F_L$ (with eigenvalues $\lambda_L \pm \ri \omega_L$)
and a stable focus for $F_R$ (with eigenvalues $\lambda_R \pm \ri \omega_R$),
then a unique limit cycle is created at the bifurcation (assuming genericity conditions are satisfied) \cite{FrPo97,SiMe07}.
The bifurcation resembles a Hopf bifurcation, but with a linear scaling law for the size of the limit cycle.
The stability of the limit cycle (and so the {\em criticality} of the bifurcation) is determined by the sign of
\begin{equation}
\alpha = \frac{\lambda_L}{\omega_L} + \frac{\lambda_R}{\omega_R}.
\label{eq:alpha}
\end{equation}
If $\alpha < 0$, the limit cycle is stable (and encircles the unstable focus);
if $\alpha > 0$, the limit cycle is unstable (and encircles the stable focus).

\begin{figure}[b!]
\begin{center}
\setlength{\unitlength}{1cm}
\begin{picture}(7.2,5.4)
\put(0,0){\includegraphics[width=7.2cm]{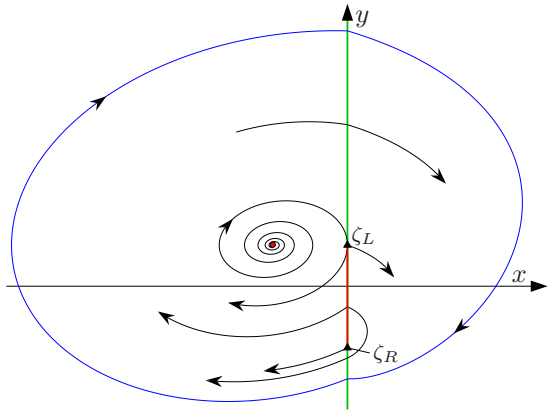}}
\put(6.72,1.69){\footnotesize $x$}
\put(4.64,5.15){\footnotesize $y$}
\put(4.88,.67){\scriptsize $\zeta_R$}
\put(4.59,2.24){\scriptsize $\zeta_L$}
\end{picture}
\caption{
A typical phase portrait of \eqref{eq:ode} subject to the assumptions described in the text.
Specifically \eqref{eq:ode} is piecewise-linear, where
$F_L$ has an unstable focus and $F_R$ has a stable focus.
Here $\alpha < 0$, the unstable focus is admissible, and there is no attracting sliding region
(instead there is a repelling sliding region with endpoints at a
visible fold $(0,\zeta_L)$ and an invisible fold $(0,\zeta_R)$).
By Theorem \ref{th:main}, there exists a unique stable limit cycle (as shown).
\label{fig:ppHLB4_4}
} 
\end{center}
\end{figure}

In this paper we show that the same result holds if \eqref{eq:ode} is discontinuous on $x=0$,
subject to an extra condition:
if $\alpha < 0$ [$\alpha > 0$] an attracting [repelling] sliding region does not coexist with the unstable [stable] focus,
see Fig.~\ref{fig:ppHLB4_4}.
To clarify, an {\em attracting sliding region} is a subset of $x=0$ where $F_L$ and $F_R$ both point towards $x=0$.
A {\em repelling sliding region} is a subset of $x=0$ where $F_L$ and $F_R$ both point away from $x=0$.
In the continuous setting, there are no sliding regions.
In the discontinuous setting, generically there is an attracting sliding region on one side of the BEB
and a repelling sliding region on the side of the BEB.
These regions shrink to a point (the boundary equilibrium) at the BEB.
The extra condition specifies which side of the BEB the sliding regions exist.
Note that in the space of two-dimensional Filippov systems, the BEB described here (termed HLB 4 in \cite{Si18c})
is not a generic codimension-one bifurcation 
because it involves two equilibria (not one equilibrium and one pseudo-equilibrium).
The merging of two foci in this fashion has been described in mathematical models
where genericity is broken by a symmetry \cite{AnVi66,ZoKu06}.

For simplicity we assume \eqref{eq:ode} is piecewise-linear.
The addition of nonlinear terms causes no qualitative change to hyperbolic equilibria
in a sufficiently small neighbourhood of the bifurcation (a simple consequence of the implicit function theorem).
The same is true for hyperbolic limit cycles, if they can be expressed as fixed points of a smooth Poincar\'e map.
For the BEB described here, this will be established formally in \cite{Si18b}.
Recently there have been many studies of two-dimensional piecewise-linear ODEs, see \cite{LlNo15} and references within.
In particular, with two foci there can exist three nested limit cycles \cite{BrMe13,FrPo14}.
Limit cycles of \eqref{eq:ode} can be analysed via a Poincar\'e map.
However, although the flow in $x<0$ and in $x>0$ is available explicitly,
Poincar\'e maps are not straight-forward to analyse when the return time of an orbit to $x=0$
cannot be obtained in closed-form, as is usually the case.
The main novelty of this paper is a generalisation of the technically difficult Poincar\'e map analysis described in \cite{FrPo98,FrPo97}
from the continuous to the discontinuous setting.

The remainder of this paper is organised as follows.
In \S\ref{sec:mainResult} we formulate the BEB in a quantitative manner and state the main result (Theorem \ref{th:main}).
We provide a minimal example and show how Theorem \ref{th:main} reduces
to the well-known continuous result in the special case that \eqref{eq:ode} is continuous on $x=0$.

In \S\ref{sec:proof} we provide a full proof of Theorem \ref{th:main}
by constructing and analysing a Poincar\'e map.
In \S\ref{sec:pseq} we determine necessary and sufficient conditions for the existence of pseudo-equilibria.
Generically there either exist no pseudo-equilibria on either side of the BEB,
or two pseudo-equilibria on both sides of the BEB.
The first case appears to occur over a wider range of parameter values.

In \S\ref{sec:threeLimitCycles} we provide an example to show that if the
extra condition described above is not satisfied then three limit cycles can be created at the BEB.
Finally \S\ref{sec:conc} provides concluding remarks.

\section{Basic properties and the main result}
\label{sec:mainResult}
\setcounter{equation}{0}

Throughout this paper we study piecewise-linear systems of the form \eqref{eq:ode}
(i.e.~$F_L$ and $F_R$ are affine functions of $x$, $y$, and $\mu$).
We refer to $(\dot{x},\dot{y}) = F_L(x,y;\mu)$ as the {\em left half-system},
and $(\dot{x},\dot{y}) = F_R(x,y;\mu)$ as the {\em right half-system}.

The BEB we wish to study involves two foci that coincide on $x=0$ at the bifurcation.
Without loss of generality we can assume the foci coincide at the origin when $\mu = 0$.
Thus $F_J(0,0;0) = (0,0)$, for each $J \in \{ L,R \}$.
Since each $F_J$ is affine, we can write
\begin{equation}
F_J(x,y;\mu) = \begin{bmatrix}
a_{1J} x + a_{2J} y + a_{3J} \mu \\
b_{1J} x + b_{2J} y + b_{3J} \mu
\end{bmatrix},
\label{eq:FJ}
\end{equation}
for some coefficients $a_{1J},\ldots,b_{3J} \in \mathbb{R}$.

The system \eqref{eq:ode} with \eqref{eq:FJ} has the property that the structure
of the dynamics is independent of the magnitude of $\mu$.
This is because if $\mu \ne 0$ then under the scaling
$(x,y) \mapsto \left( \frac{x}{|\mu|}, \frac{y}{|\mu|} \right)$
the system is unchanged except the value of $\mu$ becomes $\frac{\mu}{|\mu|} = \pm 1$.
Every bounded invariant set of \eqref{eq:ode}, such as an equilibrium or a limit cycle,
shrinks linearly to the origin as $\mu \to 0$ and
to analyse \eqref{eq:ode} it suffices to consider $\mu \in \{ -1,0,1 \}$.

Here we first clarify the BEB described above and compute equilibria, \S\ref{sub:eq}.
We then identify folds and sliding regions in \S\ref{sub:folds} (sliding motion is described in \S\ref{sec:pseq}).
We then state Theorem \ref{th:main}, \S\ref{sub:main},
illustrate the result with a simple example, \S\ref{sub:ex},
and look at the case that \eqref{eq:ode} is continuous on $x=0$, \S\ref{sub:continuous}.

\subsection{Equilibria}
\label{sub:eq}

Let
\begin{equation}
A_J = \begin{bmatrix} a_{1J} & a_{2J} \\ b_{1J} & b_{2J} \end{bmatrix}
\label{eq:AJ}
\end{equation}
denote the Jacobian matrix of \eqref{eq:FJ}.
If $\det(A_J) \ne 0$, then $(\dot{x},\dot{y}) = F_J(x,y;\mu)$ has the unique equilibrium
\begin{equation}
\begin{bmatrix} x_J^*(\mu) \\ y_J^*(\mu) \end{bmatrix} =
-A_J^{-1} \begin{bmatrix} a_{3J} \\ b_{3J} \end{bmatrix} \mu.
\label{eq:eq}
\end{equation}
The eigenvalues associated with $\left( x_J^*(\mu), y_J^*(\mu) \right)$ are those of $A_J$.
As discussed in \S\ref{sec:intro}, we suppose
\begin{equation}
\begin{split}
{\rm eig}(A_L) &= \lambda_L \pm \ri \omega_L \,, \quad
{\rm with~} \lambda_L > 0, \,\omega_L > 0, \\
{\rm eig}(A_R) &= \lambda_R \pm \ri \omega_R \,, \quad
{\rm with~} \lambda_R < 0, \,\omega_R > 0,
\end{split}
\label{eq:eigCond}
\end{equation}
so that the left half-system has an unstable focus,
and the right half-system has a stable focus.

The foci are only equilibria of \eqref{eq:ode} if they are located on the `correct' side of $x=0$.
Specifically if $x_L^*(\mu) < 0$, then $\left( x_L^*(\mu), y_L^*(\mu) \right)$
is an equilibrium of \eqref{eq:ode} and said to be {\em admissible}.
If instead $x_L^*(\mu) > 0$, then $\left( x_L^*(\mu), y_L^*(\mu) \right)$
is said to be {\em virtual}.
Similarly $\left( x_R^*(\mu), y_R^*(\mu) \right)$ is admissible if $x_R^*(\mu) > 0$,
and virtual if $x_R^*(\mu) < 0$.

\begin{figure}[t!]
\begin{center}
\setlength{\unitlength}{1cm}
\begin{picture}(7.1,11.8)
\put(0,8.2){\includegraphics[width=4.8cm]{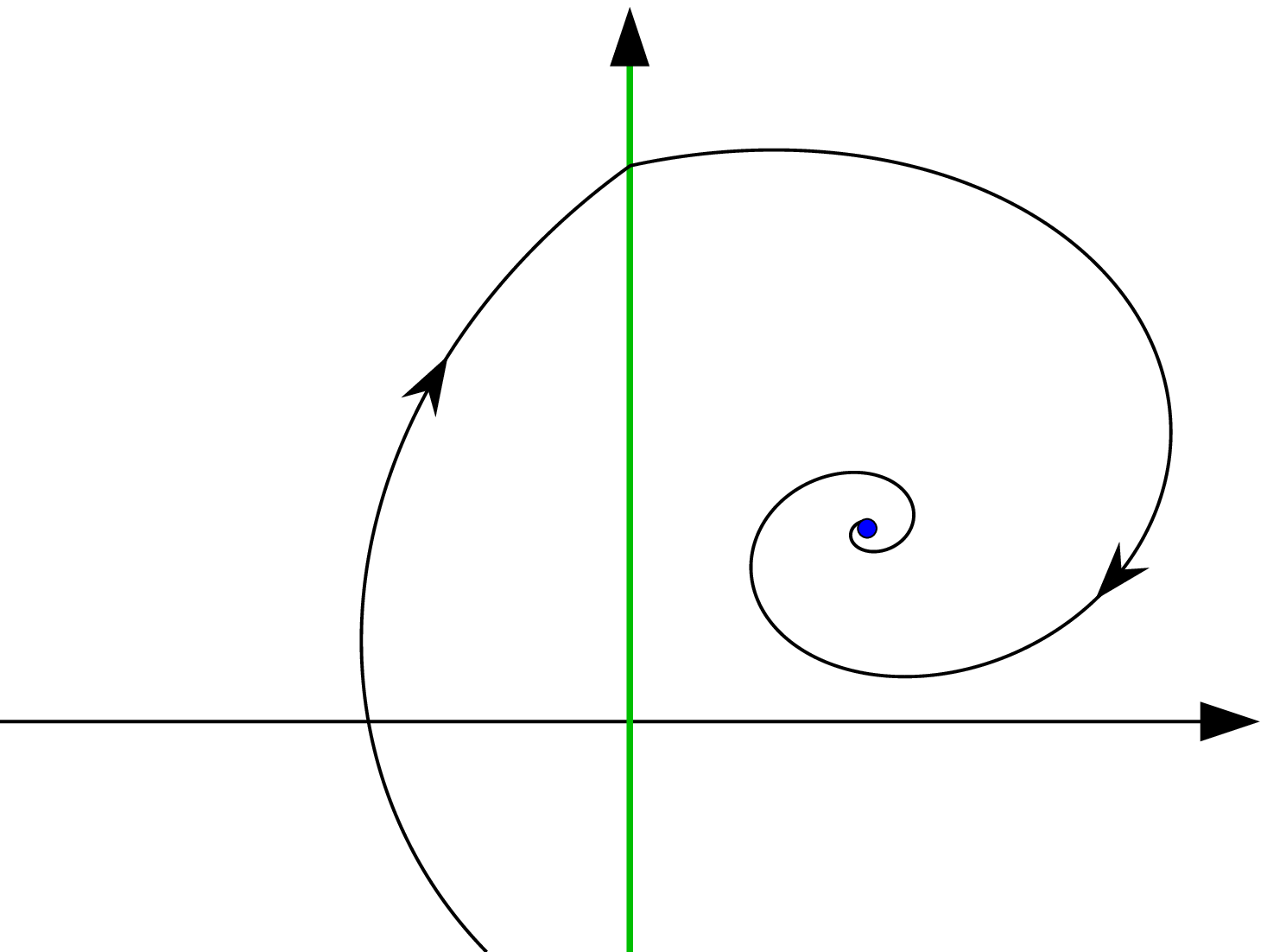}}
\put(0,4.1){\includegraphics[width=4.8cm]{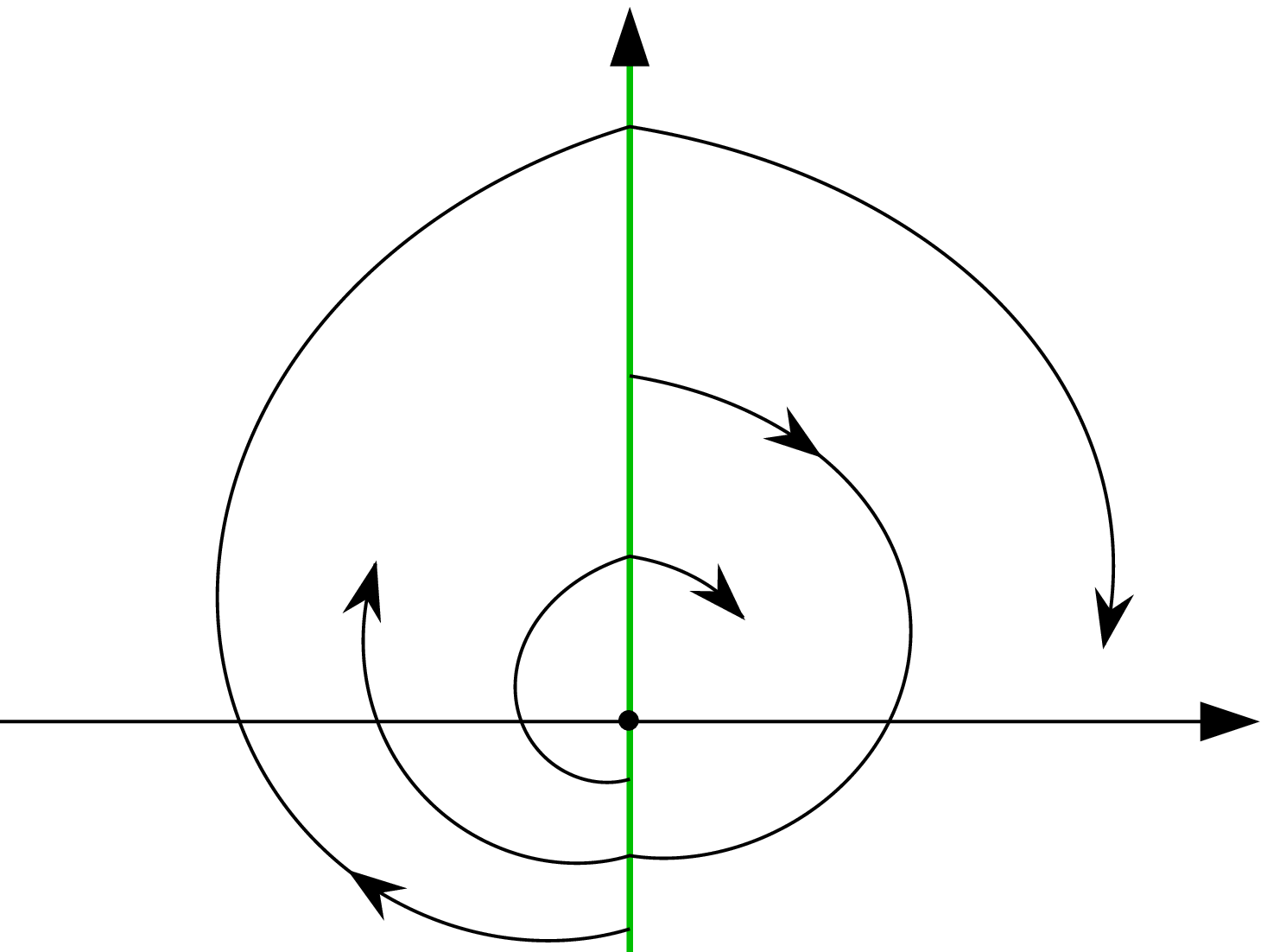}}
\put(0,0){\includegraphics[width=4.8cm]{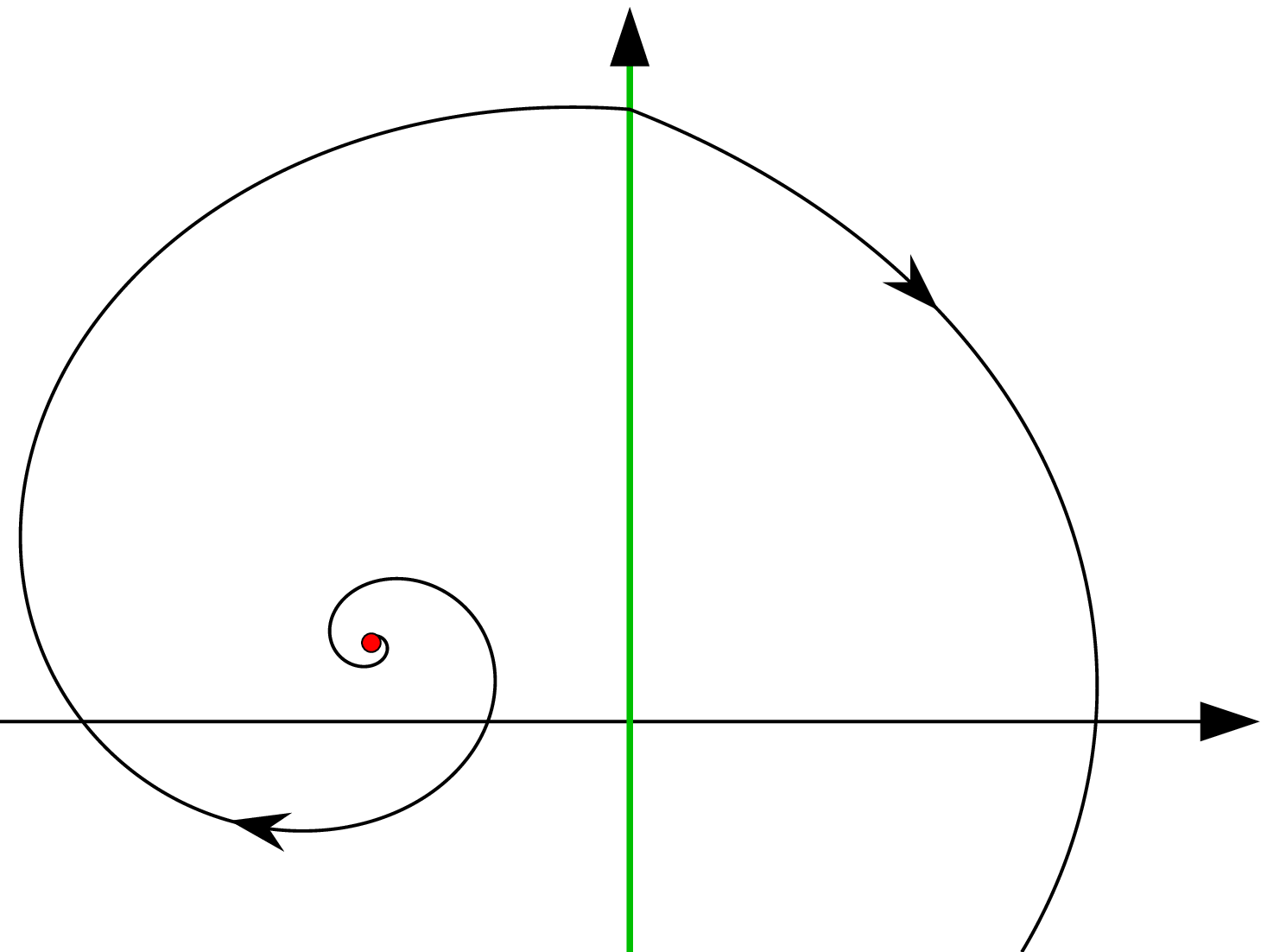}}
\put(4.32,9.14){\footnotesize $x$}
\put(2.5,11.55){\footnotesize $y$}
\put(5.3,10.08){\small {\bf a)}\hspace{2mm}$\mu < 0$}
\put(4.32,5.04){\footnotesize $x$}
\put(2.5,7.45){\footnotesize $y$}
\put(5.3,5.98){\small {\bf b)}\hspace{2mm}$\mu = 0$}
\put(4.32,.94){\footnotesize $x$}
\put(2.5,3.35){\footnotesize $y$}
\put(5.3,1.88){\small {\bf c)}\hspace{2mm}$\mu > 0$}
\end{picture}
\caption{
The basic dynamics of \eqref{eq:ode} with \eqref{eq:FJ},
assuming \eqref{eq:eigCond} is satisfied, and $a_{2L}, a_{2R}, \beta_L, \beta_R > 0$.
The right half-system has a clockwise rotating stable focus (admissible for $\mu < 0$);
the left half-system has a clockwise rotating unstable focus (admissible for $\mu > 0$).
The foci coincide at the origin when $\mu = 0$.
\label{fig:ppHLB4_6}
} 
\end{center}
\end{figure}

From \eqref{eq:eq} we obtain
$x_J^*(\mu) = \frac{-\beta_J \mu}{\lambda_J^2 + \omega_J^2}$, where
\begin{equation}
\beta_J = a_{3J} b_{2J} - a_{2J} b_{3J} \,.
\label{eq:betaJ}
\end{equation}
Therefore $\beta_J \ne 0$ ensures that $\left( x_J^*(\mu), y_J^*(\mu) \right)$
is admissible for exactly one sign of $\mu$ (i.e.~either for $\mu < 0$ or for $\mu > 0$).
This paper concerns the case that the equilibria are
admissible for different signs of $\mu$, thus $\beta_L \beta_R > 0$.
In this case \eqref{eq:ode} appears to have a single focus that changes stability
as the value of $\mu$ is varied through $0$,
much like a Hopf bifurcation.
In view of the replacement $\mu \mapsto -\mu$, we can assume $\beta_L > 0$ and $\beta_R > 0$.
Then the stable focus $\left( x_R^*(\mu), y_R^*(\mu) \right)$ is admissible for $\mu < 0$,
and the unstable focus $\left( x_L^*(\mu), y_L^*(\mu) \right)$ is admissible for $\mu > 0$, see Fig.~\ref{fig:ppHLB4_6}.

In order for a limit cycle to be created we need to
assume that the foci involve the same direction of rotation.
This assumption is equivalent to $a_{2L} a_{2R} > 0$.
If $a_{2L} > 0$ and $a_{2R} > 0$, then orbits rotate clockwise (as in Fig.~\ref{fig:ppHLB4_6}),
while if $a_{2L} < 0$ and $a_{2R} < 0$ then orbits rotate anti-clockwise.

\subsection{Folds and sliding regions}
\label{sub:folds}

Subsets of $x=0$ where $F_L$ and $F_R$ both point towards [away from] $x=0$ are called attracting [repelling] sliding regions.
Endpoints of sliding regions are usually {\em folds} where $F_L$ or $F_R$ is tangent to $x=0$.
We have assumed $a_{2L}, a_{2R} \ne 0$,
thus each half-system has a unique fold located at $(x,y) = \left( 0, \zeta_J(\mu) \right)$,
where
\begin{equation}
\zeta_J(\mu) = -\frac{a_{3J} \mu}{a_{2J}},
\label{eq:zetaJ}
\end{equation}
see Fig.~\ref{fig:ppHLB4_4}.
In this figure the fold $\left( 0, \zeta_L \right)$ is {\em visible} \cite{DiBu08}\removableFootnote{
Specifically pg.~237 where folds are called tangent points.
}
because, locally, the orbit of the left half-system that passes through $\left( 0, \zeta_L \right)$ is located in $x \le 0$
and thus is an orbit of \eqref{eq:ode}.
In contrast, $\left( 0, \zeta_R \right)$ is an {\em invisible} fold.

With $\mu = 0$ the folds coincide at the origin and there are no sliding regions.
This is because $a_{2L} a_{2R} > 0$, so the negative $y$-axis and the positive $y$-axis are both {\em crossing regions},
see Fig.~\ref{fig:ppHLB4_6}b.

The difference in the $y$-values of the folds is
$\zeta_L - \zeta_R = \frac{\gamma \mu}{a_{2L} a_{2R}}$, where 
\begin{equation}
\gamma = a_{2L} a_{3R} - a_{3L} a_{2R} \,.
\label{eq:gamma}
\end{equation}
Thus the condition $\gamma \ne 0$ ensures a sliding region exists for all $\mu \ne 0$.
It is straight-forward to show that this region is attracting for one sign of $\mu$ and repelling for the other sign of $\mu$
and we provide the following lemma without proof\removableFootnote{
Here is a proof:
In view of the above discussion it remains to classify the sliding region as attracting or repelling.
At the centre of the sliding region we have
\begin{align*}
f_L \left( 0, \frac{\zeta_L(\mu)+\zeta_R(\mu)}{2}; \mu \right) &= -\frac{\gamma \mu}{2 a_{2R}}, \\
f_R \left( 0, \frac{\zeta_L(\mu)+\zeta_R(\mu)}{2}; \mu \right) &= \frac{\gamma \mu}{2 a_{2L}},
\end{align*}
where $f_J$ denotes the first component of $F_J$.
Here
$$
-{\rm sgn}(f_L) = {\rm sgn}(f_R) = {\rm sgn} \big( a_{2L} \gamma \mu \big),
$$
and the result follows.
\manualEndProof
}.

\begin{lemma}
Consider \eqref{eq:ode} with \eqref{eq:FJ} and suppose $a_{2L} a_{2R} > 0$.
If $\gamma = 0$ or $\mu = 0$, then \eqref{eq:ode} has no sliding regions.
If $\gamma \ne 0$ and $\mu \ne 0$, then \eqref{eq:ode} has one sliding region
with endpoints at $y = \zeta_L(\mu)$ and $y = \zeta_R(\mu)$, given by \eqref{eq:zetaJ}.
The sliding region is attracting if $a_{2L} \gamma \mu < 0$, and repelling if $a_{2L} \gamma \mu > 0$.
\label{le:slidingRegion}
\end{lemma}

\subsection{A Hopf-like boundary equilibrium bifurcation}
\label{sub:main}

\begin{figure}[b!]
\begin{center}
\setlength{\unitlength}{1cm}
\begin{picture}(4.8,3.6)
\put(0,0){\includegraphics[width=4.8cm]{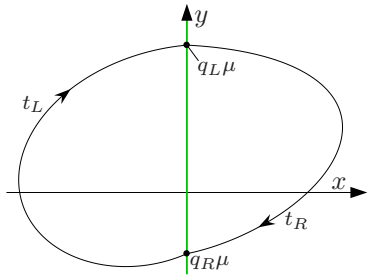}}
\put(4.32,1.13){\footnotesize $x$}
\put(2.5,3.35){\footnotesize $y$}
\put(2.54,2.7){\scriptsize $q_L \mu$}
\put(2.44,.12){\scriptsize $q_R \mu$}
\put(3.7,.65){\scriptsize $t_R$}
\put(.21,2.2){\scriptsize $t_L$}
\end{picture}
\caption{
The limit cycle of Theorem \ref{th:main} in the case $a_{2L}, a_{2R} > 0$ (clockwise rotation).
The limit cycle intersects $x=0$ at $y = q_L \mu$ and $y = q_R \mu$.
Its evolution time in $x<0$ is $t_L$, and in $x>0$ is $t_R$.
These values are given implicitly by \eqref{eq:qpmtpmImplicit}.
\label{fig:ppHLB4_5}
} 
\end{center}
\end{figure}

Here we state our main result for the existence of a unique limit cycle.
We provide formulas for its evolution time in $x<0$ (denoted $t_L$),
and in $x>0$ (denoted $t_R$),
and its points of intersection with $x=0$.
Since \eqref{eq:ode} is piecewise-linear,
$t_L$ and $t_R$ are independent of $\mu$ and
the size of the limit cycle is proportional to $|\mu|$.
There exist $q_L, q_R \in \mathbb{R}$ such that as time increases
the limit cycle crosses from $x<0$ to $x>0$ at $(0,q_L \mu)$, 
and crosses from $x>0$ to $x<0$ at $(0,q_R \mu)$, see Fig.~\ref{fig:ppHLB4_5}.

\begin{theorem}
Consider \eqref{eq:ode} with \eqref{eq:FJ}.
Suppose \eqref{eq:eigCond} is satisfied,
$\beta_L > 0$, $\beta_R > 0$, $a_{2L} a_{2R} > 0$, and $a_{2L} \gamma \ge 0$.
Then \eqref{eq:ode} has
\begin{enumerate}
\item
a stable focus in $x>0$ for $\mu < 0$,
an unstable focus in $x<0$ for $\mu > 0$, and
\item
if $\alpha < 0$ [$\alpha > 0$] there exists a unique stable [unstable] limit cycle for $\mu > 0$ [$\mu < 0$],
and no limit cycle for $\mu < 0$ [$\mu > 0$].
\end{enumerate}
The quantities $q_L, q_R, t_L, t_R \in \mathbb{R}$ satisfy\removableFootnote{
In the proof these are derived assuming $a_{2L} > 0$ and $\mu = 1$,
but they hold for any $a_{2L} \ne 0$ and $\mu \ne 0$, see {\sc checkFormulasHLB4.m}.
The trick is to remember that the actual $y$-values are $y_R = q_R \mu$ and $y_L = q_L \mu$.
}
\begin{equation}
\begin{split}
q_R + \frac{a_{3R}}{a_{2R}} &= \frac{\xi_R \,\re^{-\lambda_R t_R}
\varrho \big( \omega_R t_R; \frac{\lambda_R}{\omega_R} \big)}{\sin(\omega_R t_R)}, \\
q_L + \frac{a_{3R}}{a_{2R}} &= \frac{-\xi_R \,\re^{\lambda_R t_R}
\varrho \big( \omega_R t_R; -\frac{\lambda_R}{\omega_R} \big)}{\sin(\omega_R t_R)}, \\
q_L + \frac{a_{3L}}{a_{2L}} &= \frac{\xi_L \,\re^{-\lambda_L t_L}
\varrho \big( \omega_L t_L; \frac{\lambda_L}{\omega_L} \big)}{\sin(\omega_L t_L)}, \\
q_R + \frac{a_{3L}}{a_{2L}} &= \frac{-\xi_L \,\re^{\lambda_L t_L}
\varrho \big( \omega_L t_L; -\frac{\lambda_L}{\omega_L} \big)}{\sin(\omega_L t_L)},
\end{split}
\label{eq:qpmtpmImplicit}
\end{equation}
where $\xi_J = \frac{\beta_J \omega_J}{a_{2J} \left( \lambda_J^2 + \omega_J^2 \right)}$\removableFootnote{
In my HLB survey, I define $\kappa_R = -\xi_R$, and $\kappa_L = \xi_L$.
},
for each $J \in \{ L, R \}$, and
\begin{equation}
\varrho(s;\nu) = 1 - \re^{\nu s} \left( \cos(s) - \nu \sin(s) \right).
\label{eq:auxFunc}
\end{equation}
\label{th:main}
\end{theorem}

Theorem \ref{th:main} is proved in \S\ref{sec:proof}.
Notice we cannot provide explicit formulas for $t_L$ and $t_R$ in terms of the parameters of \eqref{eq:ode}.
Instead they are given implicitly by \eqref{eq:qpmtpmImplicit}
in terms of the auxiliary function $\varrho(s;\nu)$ which was introduced in \cite{AnVi66}.

\subsection{An example}
\label{sub:ex}

As a simple example consider the system
\begin{equation}
\begin{bmatrix} \dot{x} \\ \dot{y} \end{bmatrix} =
\begin{cases}
\begin{bmatrix}
y \\
-x + 2 \lambda_L y - \mu
\end{bmatrix}, & x < 0, \\
\begin{bmatrix}
-x + y + \mu \\
-x - \mu
\end{bmatrix}, & x > 0.
\end{cases}
\label{eq:ex1}
\end{equation}
This is of the form \eqref{eq:ode} with \eqref{eq:FJ}.
The right half-system has a stable focus,
and the left half-system has an unstable focus when $0 < \lambda_L < 1$.

\begin{figure}[t!]
\begin{center}
\setlength{\unitlength}{1cm}
\begin{picture}(9.5,17.2)
\put(0,11.8){\includegraphics[width=7.2cm]{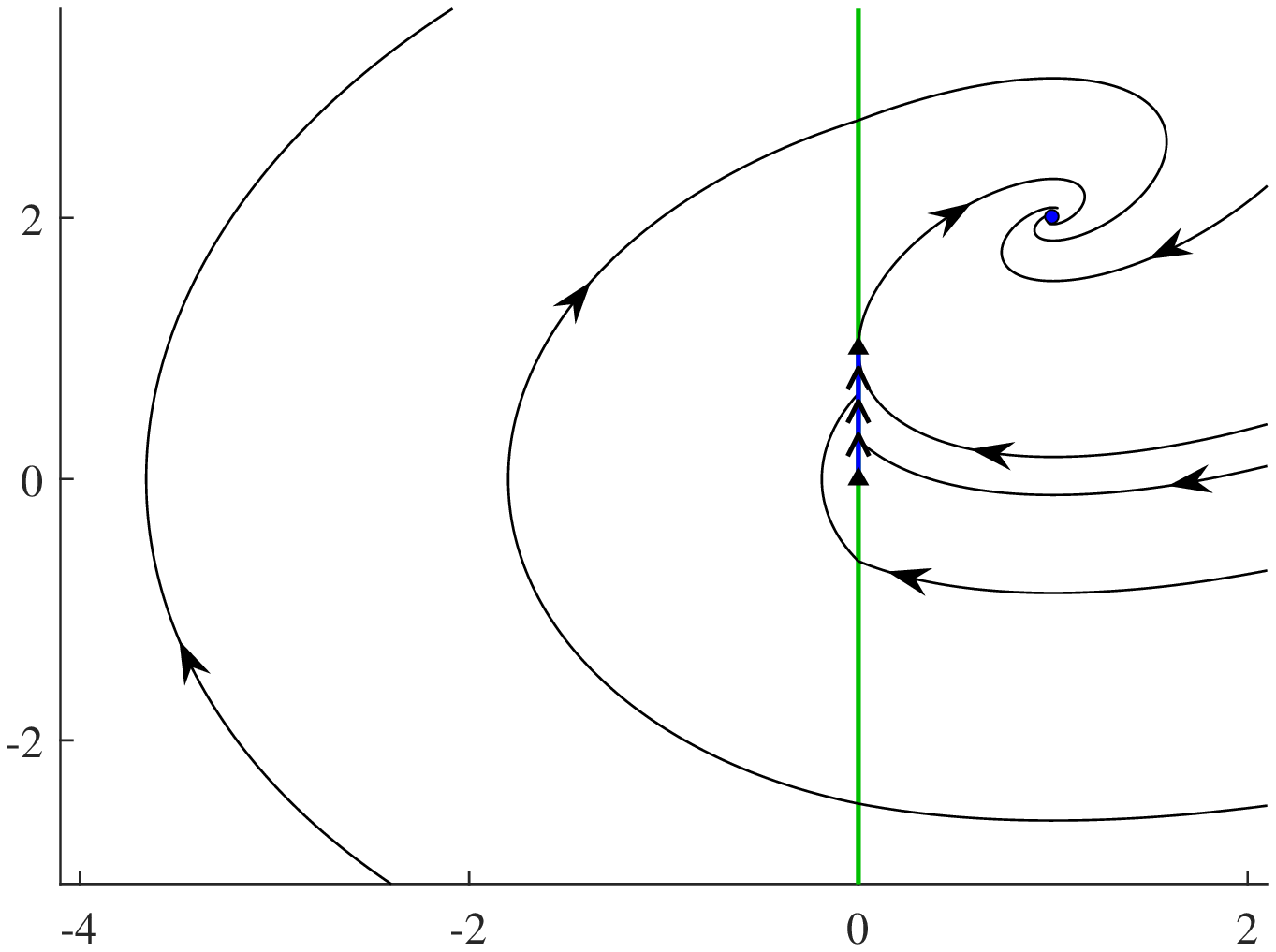}}
\put(0,5.9){\includegraphics[width=7.2cm]{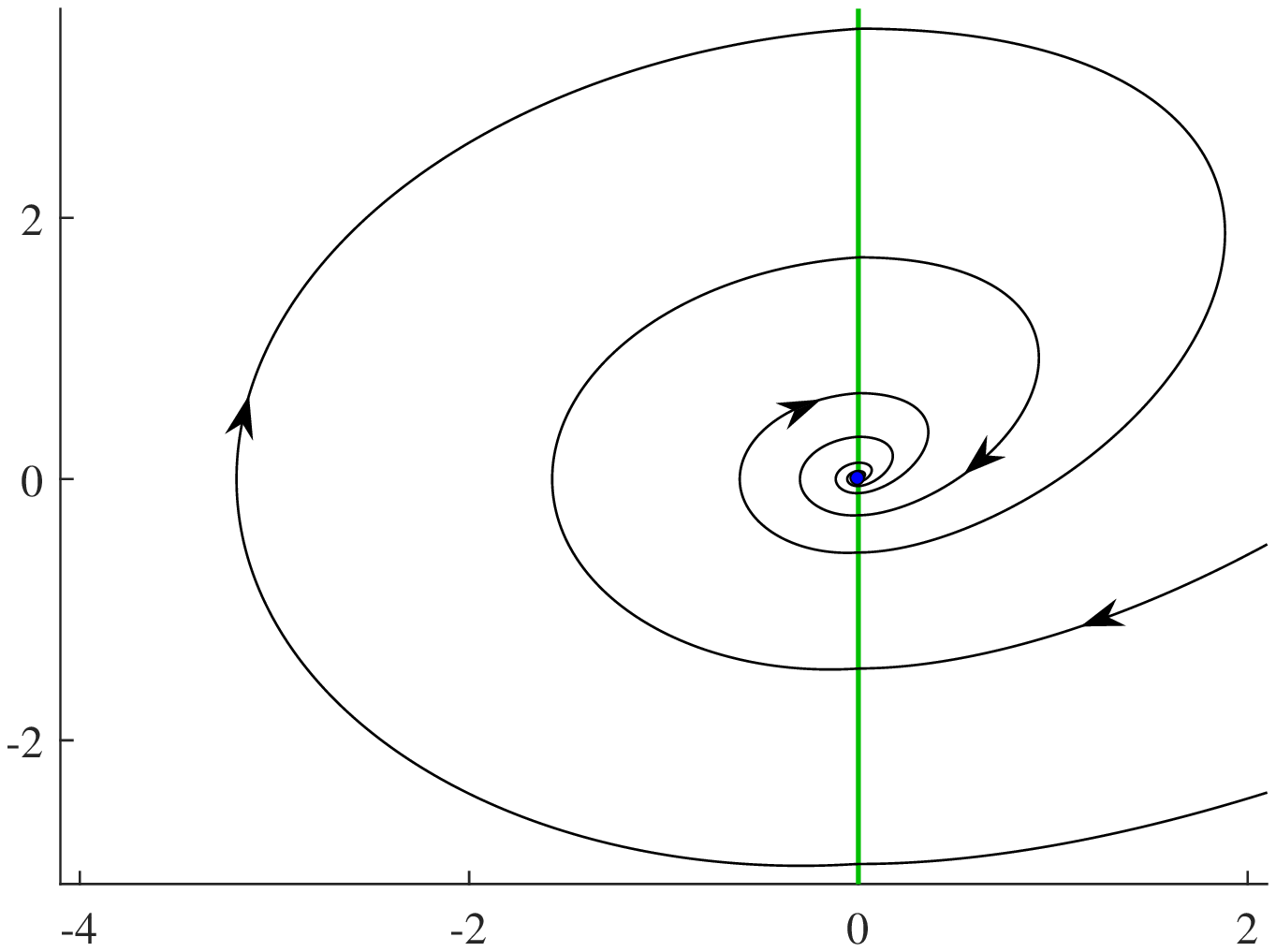}}
\put(0,0){\includegraphics[width=7.2cm]{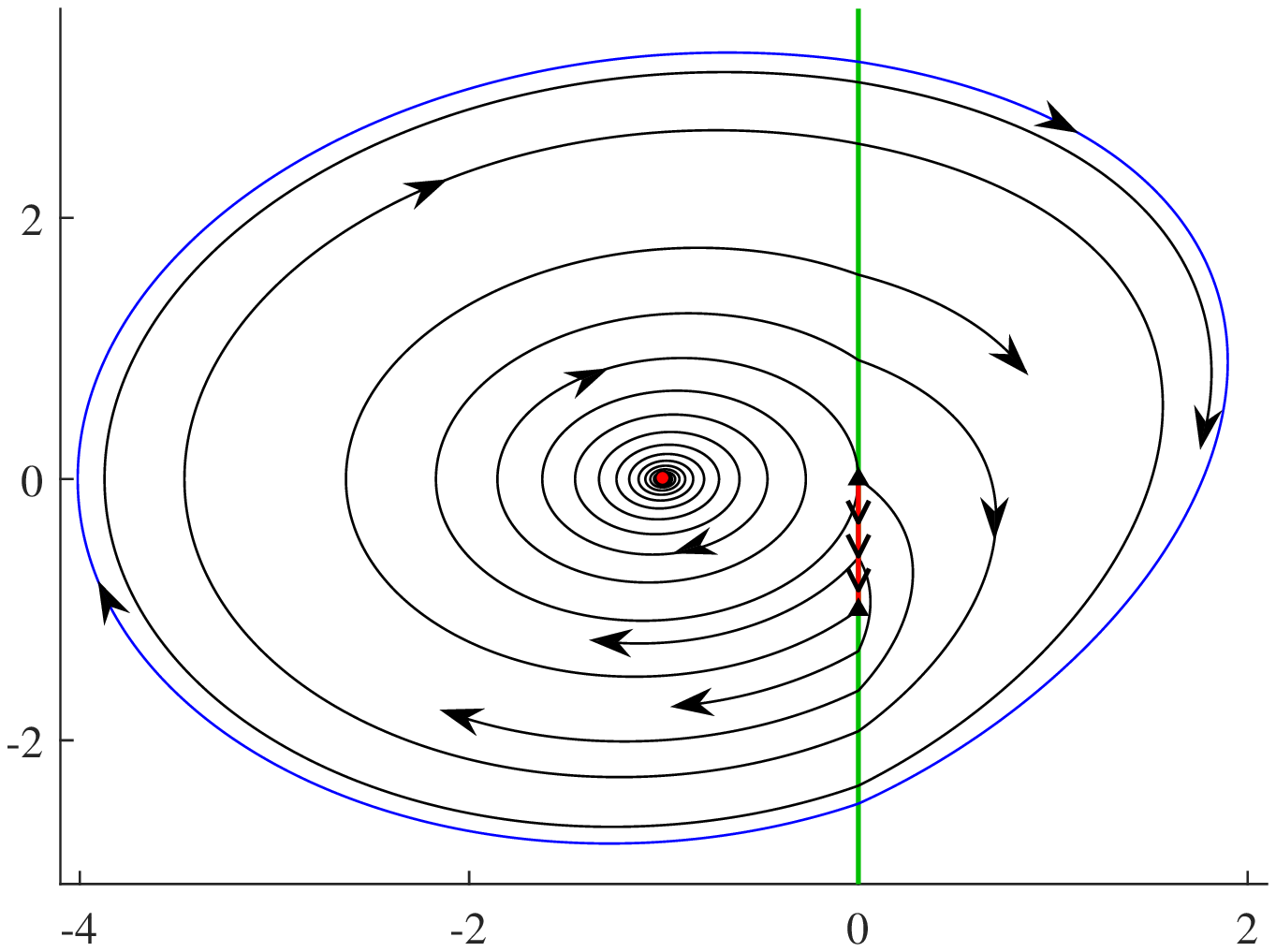}}
\put(4.87,11.8){\small $x$}
\put(0,14.62){\small $y$}
\put(7.7,14.62){\small {\bf a)}\hspace{2mm}$\mu = -1$}
\put(4.87,5.9){\small $x$}
\put(0,8.72){\small $y$}
\put(7.7,8.72){\small {\bf b)}\hspace{2mm}$\mu = 0$}
\put(4.87,0){\small $x$}
\put(0,2.82){\small $y$}
\put(7.7,2.82){\small {\bf c)}\hspace{2mm}$\mu = 1$}
\end{picture}
\caption{
Phase portraits of \eqref{eq:ex1} with $\lambda_L = 0.05$.
When $\mu < 0$ there exists a stable focus and an attracting sliding region.
When $\mu > 0$ there exists an unstable focus, a repelling sliding region, and a stable limit cycle.
\label{fig:ppHLB4_1}
} 
\end{center}
\end{figure}

Here $a_{2L} = a_{2R} = \beta_L = \beta_R = \gamma = 1$,
thus \eqref{eq:ex1} satisfies the conditions of Theorem \ref{th:main}.
Fig.~\ref{fig:ppHLB4_1} shows phase portraits using $\lambda_L = 0.05$.
Here $\alpha < 0$ and so, by Theorem \ref{th:main}, a unique stable limit cycle exists for $\mu > 0$.

\subsection{Continuous piecewise-smooth systems}
\label{sub:continuous}

Here we suppose \eqref{eq:ode} is continuous on $x=0$.
That is, $F_L(0,y;\mu) = F_R(0,y;\mu)$, for all $y$ and $\mu$.
Given that $F_L$ and $F_R$ take the form \eqref{eq:FJ},
we must have $a_{2L} = a_{2R}$, $a_{3L} = a_{3R}$, $b_{2L} = b_{2R}$, and $b_{3L} = b_{3R}$.

Again suppose that the eigenvalues of $A_L$ and $A_R$ satisfy \eqref{eq:eigCond}
so that the left half-system has an unstable focus and the right half-system has a stable focus.
The foci are located at the origin when $\mu = 0$.
Assuming the unstable focus moves away from $x=0$ as the value of $\mu$ is varied from $0$,
then $\beta_L \ne 0$ and we can assume $\beta_L > 0$ if we allow the replacement $\mu \mapsto -\mu$.

In this situation all conditions of Theorem \ref{th:main} are satisfied by continuity.
Specifically, by \eqref{eq:eigCond} we have $a_{2L} \ne 0$,
so since $a_{2L} = a_{2R}$ we have $a_{2L} a_{2R} > 0$ (i.e.~the foci have the same direction of rotation).
By continuity, $\beta_R = \beta_L > 0$.
Also, by continuity, $\gamma = 0$.
Thus a unique limit cycle is created and its stability is determined by the sign of $\alpha$.
Therefore, in the special case that \eqref{eq:ode} is continuous,
Theorem \ref{th:main} reduces to the Hopf-like bifurcation theorem of \cite{FrPo97,SiMe07} for continuous systems.

\section{Proof of Theorem \ref{th:main}}
\label{sec:proof}
\setcounter{equation}{0}

\myStep{1}{Sign assumptions and the Poincar\'e map $P = P_L \circ P_R$.}\\
By symmetry it suffices to prove the result for $\mu > 0$.
This is justified through the change of variables
$(x,y;\mu;t) \mapsto (-x,y;-\mu;-t)$ which flips the sign of $\mu$ and
transforms \eqref{eq:ode} into another system
satisfying the conditions of Theorem \ref{th:main}.
In fact it suffices to consider $\mu = 1$
in view of scaling property discussed at the start of \S\ref{sec:mainResult},
and so for the remainder of the proof we assume $\mu = 1$.

By assumption $a_{2L} \ne 0$.
Without loss of generality we may assume $a_{2L} > 0$
(justified by the change of variables $y \mapsto -y$, which also flips the sign of $\gamma$).
Then, by assumption, $a_{2R} > 0$ and $\gamma \ge 0$.

Since $\mu = 1$, the left half-system has a visible fold at $y = \zeta_L = -\frac{a_{3L}}{a_{2L}}$
and the right half-system has an invisible fold at $y = \zeta_R = -\frac{a_{3R}}{a_{2R}}$,
as in Fig.~\ref{fig:ppHLB4_4}.
Notice $\zeta_L - \zeta_R = \frac{\gamma}{a_{2L} a_{2R}} \ge 0$.

\begin{figure}[t!]
\begin{center}
\setlength{\unitlength}{1cm}
\begin{picture}(7.2,5.4)
\put(0,0){\includegraphics[width=7.2cm]{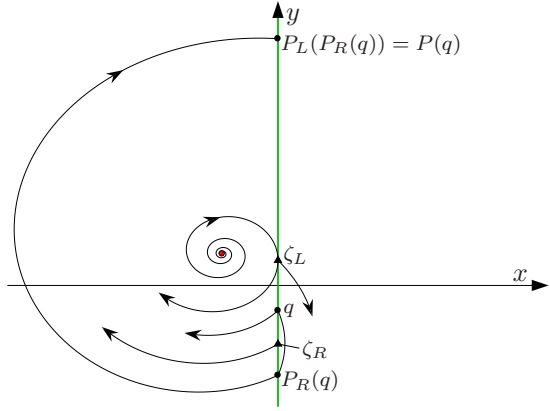}}
\put(6.72,1.67){\footnotesize $x$}
\put(3.71,5.15){\footnotesize $y$}
\put(3.64,.27){\scriptsize $P_R(q)$}
\put(3.92,.69){\scriptsize $\zeta_R$}
\put(3.67,1.25){\scriptsize $q$}
\put(3.67,1.95){\scriptsize $\zeta_L$}
\put(3.64,4.75){\scriptsize $P_L(P_R(q)) = P(q)$}
\end{picture}
\caption{
A phase portrait of \eqref{eq:ode}, subject to the conditions of Theorem \ref{th:main},
illustrating the Poincar\'e map $P(q)$.
\label{fig:schemHLB4}
} 
\end{center}
\end{figure}

Given $q > \zeta_R$, consider the forward orbit of $(x,y) = (0,q)$
that immediately enters $x > 0$, see Fig.~\ref{fig:schemHLB4}.
Let $P_R(q)$ denote the $y$-value of the next intersection of this orbit with $x = 0$,
and let $T_R(q)$ denote the corresponding evolution time.
Similarly given $q < \zeta_L$, consider the forward orbit of $(x,y) = (0,q)$
that immediately enters $x < 0$.
Let $P_L(q)$ denote the $y$-value of the next intersection of this orbit with $x = 0$,
and let $T_L(q)$ denote the corresponding evolution time.
Notice $P_R(q) < \zeta_R \le \zeta_L$, thus the Poincar\'e map
\begin{equation}
P(q) = P_L \left( P_R(q) \right),
\label{eq:P}
\end{equation}
is well-defined for all $q > \zeta_R$.

\myStep{2}{Formulas for $P_R$ and $T_R$.}\\
The right half-system has the unique equilibrium
\begin{equation}
\begin{bmatrix} x_R^* \\ y_R^* \end{bmatrix} = -A_R^{-1}
\begin{bmatrix} a_{3R} \\ b_{3R} \end{bmatrix}.
\label{eq:zRstar}
\end{equation}
The flow of the right half-system is
\begin{equation}
\begin{bmatrix} \varphi^R_t(x,y) \\ \psi^R_t(x,y) \end{bmatrix} =
\re^{t A_R} \left( \begin{bmatrix} x \\ y \end{bmatrix}
- \begin{bmatrix} x_R^* \\ y_R^* \end{bmatrix} \right)
+ \begin{bmatrix} x_R^* \\ y_R^* \end{bmatrix},
\label{eq:phipsi}
\end{equation}
where
\begin{equation}
\re^{t A_R} = \re^{\lambda_R t} \begin{bmatrix}
\cos(\omega_R t) + \frac{a_{1R}-b_{2R}}{2 \omega_R} \,\sin(\omega_R t) &
\frac{a_{2R}}{\omega_R} \,\sin(\omega_R t) \\
\frac{b_{1R}}{\omega_R} \,\sin(\omega_R t) &
\cos(\omega_R t) - \frac{a_{1R}-b_{2R}}{2 \omega_R} \,\sin(\omega_R t) \end{bmatrix},
\label{eq:exptAR}
\end{equation}
and
\begin{equation}
\begin{split}
\lambda_R &= \frac{a_{1R} + b_{2R}}{2}, \\
\omega_R &= \sqrt{-a_{2R} b_{1R} - \frac{(a_{1R}-b_{2R})^2}{4}}.
\end{split}
\nonumber
\end{equation}
Upon substituting $(x,y) = (0,q)$ into \eqref{eq:phipsi},
we obtain, after much simplification,
\begin{align}
\varphi^R_t(0,q) &=
\frac{a_{2R}}{\omega_R} \,\re^{\lambda_R t} \sin(\omega_R t) \left( q - \zeta_R \right)
- \frac{a_{2R} \xi_R}{\omega_R} \,\varrho \left( \omega_R t; \frac{\lambda_R}{\omega_R} \right),
\label{eq:phi2R} \\
\psi^R_t(0,q) &=
\left( \cos(\omega_R t) - \frac{a_{1R}-b_{2R}}{2 \omega_R} \,\sin(\omega_R t) \right)
\left( \re^{\lambda_R t} \left( q - \zeta_R \right) - \frac{\xi_R \varrho
\left( \omega_R t; \frac{\lambda_R}{\omega_R} \right)}{\sin(\omega_R t)} \right) \nonumber \\
&\quad- \frac{\xi_R \re^{\lambda_R t} \varrho \left( \omega_R t; -\frac{\lambda_R}{\omega_R} \right)}
{\sin(\omega_R t)} + \zeta_R \,,
\label{eq:psi2R}
\end{align}
where $\xi_R$ and $\varrho(s;\nu)$ are defined in the theorem statement.
By definition, $\varphi^R_{T_R(q)}(0,q) = 0$, thus by \eqref{eq:phi2R} we have
\begin{equation}
q - \zeta_R = \frac{\xi_R \,\re^{-\lambda_R T_R}
\varrho \left( \omega_R T_R; \frac{\lambda_R}{\omega_R} \right)}
{\sin(\omega_R T_R)}.
\label{eq:TRimplicit}
\end{equation}
Also $\psi^R_{T_R(q)}(0,q) = P_R(q)$,
thus from the way we have factored \eqref{eq:psi2R} we immediately obtain
\begin{equation}
P_R - \zeta_R = \frac{-\xi_R \,\re^{\lambda_R T_R}
\varrho \left( \omega_R T_R; -\frac{\lambda_R}{\omega_R} \right)}
{\sin(\omega_R T_R)}.
\label{eq:PRimplicit}
\end{equation}

\myStep{3}{Derivatives of $P_R$ and $T_R$.}\\
By using the identity
\begin{equation}
\frac{\partial}{\partial s} \frac{\re^{-\nu s} \varrho(s;\nu)}{\sin(s)}
= \frac{\varrho(s;-\nu)}{\sin^2(s)},
\label{eq:auxFuncFuncDeriv}
\end{equation}
to differentiate \eqref{eq:TRimplicit}, we obtain
\begin{equation}
\frac{d T_R}{d q} = -\frac{\re^{\lambda_R T_R} \sin(\omega_R T_R)}
{\omega_R \left( P_R - \zeta_R \right)}.
\label{eq:dTRdq}
\end{equation}
Since $\left( x_R^*, y_R^* \right)$ is virtual,
the orbit of the right half-system from $(0,q)$ to $(0,P_R(q))$ completes less than half a revolution
about $\left( x_R^*, y_R^* \right)$, hence $T_R(q) \in \left( 0, \frac{\pi}{\omega_R} \right)$.
Also $\omega_R > 0$ and $P_R(q) < \zeta_R$, thus,
by \eqref{eq:dTRdq}, $T_R(q)$ is an increasing function of $q$.
It follows that $T_R(q) \to \frac{\pi}{\omega_R}$ as $q \to \infty$ by \eqref{eq:TRimplicit}.

By \eqref{eq:TRimplicit} and \eqref{eq:PRimplicit} we have
\begin{equation}
\frac{P_R - \zeta_R}{q - \zeta_R} = \frac{-\re^{2 \lambda_R T_R} \varrho \left( \omega_R T_R; -\frac{\lambda_R}{\omega_R} \right)}
{\varrho \left( \omega_R T_R; \frac{\lambda_R}{\omega_R} \right)}.
\nonumber
\end{equation}
Substituting $T_R = \frac{\pi}{\omega_R}$ gives $\frac{P_R - \zeta_R}{q - \zeta_R} = -\re^{\frac{\lambda_R \pi}{\omega_R}}$,
by using the definition of $\varrho(s;\nu)$.
Thus $P_R(q) \sim -q \,\re^{\frac{\lambda_R \pi}{\omega_R}}$ as $q \to \infty$.
By applying \eqref{eq:auxFuncFuncDeriv} to both \eqref{eq:TRimplicit} and \eqref{eq:PRimplicit} we obtain
\begin{equation}
\frac{d P_R}{d q} = -\frac{\varrho \left( \omega_R T_R; \frac{\lambda_R}{\omega_R} \right)}
{\varrho \left( \omega_R T_R; -\frac{\lambda_R}{\omega_R} \right)},
\label{eq:dPRdq}
\end{equation}
and a further application of \eqref{eq:TRimplicit} and \eqref{eq:PRimplicit} produces
\begin{equation}
\frac{d P_R}{d q} = \frac{q - \zeta_R}{P_R - \zeta_R} \,\re^{2 \lambda_R T_R}.
\label{eq:dPRdq2}
\end{equation}
Since $P_R(q) < \zeta_R < q$, we conclude that $P_R(q)$ is a decreasing function of $q$.

\myStep{4}{Formulas for $P_L$ and $T_L$ and their derivatives.}\\
By repeating the above analysis for the left half-system we obtain
\begin{equation}
q - \zeta_L = \frac{\xi_L \,\re^{-\lambda_L T_L}
\varrho \left( \omega_L T_L; \frac{\lambda_L}{\omega_L} \right)}
{\sin(\omega_L T_L)}.
\label{eq:TLimplicit}
\end{equation}
and
\begin{equation}
P_L - \zeta_L = \frac{-\xi_L \,\re^{\lambda_L T_L}
\varrho \left( \omega_L T_L; -\frac{\lambda_L}{\omega_L} \right)}
{\sin(\omega_L T_L)}.
\label{eq:PLimplicit}
\end{equation}
Since $\left( x_L^*, y_L^* \right)$ is admissible,
the orbit of the left half-system from $(0,q)$ to $(0,P_L(q))$ completes more than half a revolution
about $\left( x_L^*, y_L^* \right)$, hence $T_L(q) \in \left( \frac{\pi}{\omega_L}, \frac{2 \pi}{\omega_L} \right)$.
It follows that $T_L(q)$ is an increasing function of $q$
with $T_L(q) \to \frac{\pi}{\omega_L}$ as $q \to -\infty$.
Also
\begin{equation}
\frac{d P_L}{d q} = \frac{q - \zeta_L}{P_L - \zeta_L} \,\re^{2 \lambda_L T_L},
\label{eq:dPLdq2}
\end{equation}
and so $P_L(q)$ is a decreasing function of $q$ with
$P_L \sim -q \re^{\frac{\lambda_L \pi}{\omega_L}}$ as $q \to -\infty$.

\myStep{5}{Properties of $P$.}\\
From the limiting values of $P_R$ and $P_L$,
we obtain $P(q) \sim q \re^{\alpha \pi}$, as $q \to \infty$,
where $\alpha = \frac{\lambda_L}{\omega_L} + \frac{\lambda_R}{\omega_R}$.
From \eqref{eq:dPRdq2} and \eqref{eq:dPLdq2}, we obtain
\begin{equation}
\frac{d P}{d q} =
\frac{\left( q - \zeta_R \right) \left( P_R - \zeta_L \right)}
{\left( P_R - \zeta_R \right) \left( P - \zeta_L \right)} \,\re^{2 h},
\label{eq:dPdq}
\end{equation}
where
\begin{equation}
h(q) = \lambda_R T_R(q) + \lambda_L T_L \left( P_R(q) \right).
\label{eq:h}
\end{equation}
Notice $h(q)$ is a decreasing function of $q$.
This is because $\lambda_R < 0$ and $T_R(q)$ is increasing,
thus the first term in \eqref{eq:h} is decreasing.
Also $\lambda_L > 0$, $T_L(q)$ is increasing, and $P_R(q)$ is decreasing, 
thus the second term in \eqref{eq:h} is also decreasing.

\myStep{6}{Demonstration that the smallest fixed point of $P$ is asymptotically stable.}\\
Suppose for a moment that $P(q)$ has a fixed point.
Let $q^*$ be the smallest such point.
Since $\lim_{q \to \zeta_R} P(q) > \zeta_L \ge \zeta_R$,
we must have $P(q) > q$ for all $\zeta_R < q < q^*$,
and thus $\frac{d P}{d q}(q^*) \le 1$.

Suppose $\frac{d P}{d q}(q^*) = 1$, for a contradiction.
Then $\frac{d^2 P}{d q^2}(q^*) \ge 0$,
but by differentiating \eqref{eq:dPdq} we obtain
\begin{equation}
\frac{d^2 P}{d q^2} = \left(
\frac{1}{q - \zeta_R}
- \frac{\frac{d P}{d q}}{P - \zeta_L}
+ \frac{\frac{d P_R}{d q}}{P_R - \zeta_L}
- \frac{\frac{d P_R}{d q}}{P_R - \zeta_R}
+ 2 \frac{d h}{d q} \right) \frac{d P}{d q},
\nonumber
\end{equation}
and upon substituting $q = P = q^*$ and
$\frac{d P}{d q} = 1$ we get
\begin{equation}
\frac{d^2 P}{d q^2}(q^*) =
\frac{-(\zeta_L - \zeta_R)}{(q^*-\zeta_R)(q^*-\zeta_L)} +
\frac{(\zeta_L - \zeta_R) \frac{d P_R}{d q}}
{(P_R-\zeta_R)(P_R-\zeta_L)} +
2 \frac{d h}{d q}.
\label{eq:secondDeriv}
\end{equation}
Since $P_R < \zeta_R \le \zeta_L < q^*$, $\frac{d P_R}{d q} < 0$, and $\frac{d h}{d q} < 0$,
the first two terms in \eqref{eq:secondDeriv} are zero (if $\gamma = 0$) or negative (if $\gamma > 0$),
and the third term is negative.
Thus $\frac{d^2 P}{d q^2}(q^*) < 0$,
which is a contradiction.
Therefore $\frac{d P}{d q}(q^*) < 1$,
and so $q^*$ is an asymptotically stable fixed point of $P(q)$.

\myStep{7}{Demonstration that fixed points of $P$ are unique.}\\
Now suppose for a contradiction that $P(q)$ has other fixed points.
Let $q^{**}$ be the next smallest fixed point.
Then $P(q) < q$ for all $q^* < q < q^{**}$
and thus $\frac{d P}{d q}(q^{**}) \ge 1$.
At a fixed point, \eqref{eq:dPdq} is satisfied with $P = q$, that is
$\frac{d P}{d q}(q) = G(q)$, where
\begin{equation}
G = \frac{(q - \zeta_R)(P_R - \zeta_L)}
{(P_R - \zeta_R)(q - \zeta_L)} \,\re^{2 h}.
\label{eq:G}
\end{equation}
By differentiating \eqref{eq:G} we obtain
\begin{equation}
\frac{d G}{d q} = \left( \frac{-(\zeta_L - \zeta_R)}{(q - \zeta_R)(q - \zeta_L)}
+ \frac{(\zeta_L - \zeta_R) \frac{d P_R}{d q}}{(P_R - \zeta_R)(P_R - \zeta_L)} + 2 \frac{d h}{d q} \right) G,
\nonumber
\end{equation}
and, analogous to the previous step (also observing $G(q) > 0$), we conclude that $\frac{d G}{d q} < 0$.
Thus the value of $\frac{d P}{d q}$ at fixed points decreases with $q$.
But $\frac{d P}{d q}(q^*) < 1$, so we cannot have $\frac{d P}{d q}(q^{**}) \ge 1$.
This is contradiction, hence $P(q)$ has no other fixed points.

\myStep{8}{Final remarks.}\\
In summary we have shown that if $P$ has a fixed point, then it is unique and asymptotically stable.
Since $\lim_{q \to \zeta_R} P(q) > \zeta_R$ and $P(q) \sim q \re^{\alpha \pi}$ as $q \to \infty$,
this cannot occur if $\alpha > 0$.
Thus if $\alpha > 0$, $P$ has no fixed points and so \eqref{eq:ode} has no limit cycles.
If $\alpha < 0$, then $P$ has a fixed point by the intermediate value theorem,
and so \eqref{eq:ode} has a unique stable limit cycle.
The fixed point of $P$ is the value $q_R$ in the theorem.
Also $q_L = P_R(q_R)$, $t_R = T_R(q_R)$, and $t_L = T_L(q_L)$.
Consequently the equations \eqref{eq:qpmtpmImplicit} follow
immediately from \eqref{eq:TRimplicit}, \eqref{eq:PRimplicit}, \eqref{eq:TLimplicit}, and \eqref{eq:PLimplicit}.
\manualEndProof

\section{Sliding motion and pseudo-equilibria}
\label{sec:pseq}
\setcounter{equation}{0}

Here we write \eqref{eq:FJ} as
\begin{equation}
F_J(x,y;\mu) =
\begin{bmatrix} f_J(x,y;\mu) \\ g_J(x,y;\mu) \end{bmatrix},
\label{eq:FJ2}
\end{equation}
for each $J \in \{ L, R \}$.
Attracting sliding regions are subsets of $x=0$ for which $f_L > 0$ and $f_R < 0$.
Repelling sliding regions are subsets of $x=0$ for which $f_L < 0$ and $f_R > 0$.
Recall from Lemma \ref{le:slidingRegion} that, assuming $\gamma \ne 0$,
\eqref{eq:ode} has one sliding region for all $\mu \ne 0$
with endpoints at $y = \zeta_L$ and $y = \zeta_R$.

\subsection{Sliding motion}

On sliding regions, sliding motion is defined most simply by constructing a sliding vector field.
Following the usual Filippov convention \cite{DiBu08,Fi88},
this vector field is the convex combination of $F_L$ and $F_R$ that is tangent to $x=0$.
We write
\begin{equation}
\begin{bmatrix} 0 \\ g_{\rm slide}(y;\mu) \end{bmatrix} =
\big( 1 - \theta(y;\mu) \big) \begin{bmatrix} f_L(0,y;\mu) \\ g_L(0,y;\mu) \end{bmatrix} +
\theta(y;\mu) \begin{bmatrix} f_R(0,y;\mu) \\ g_R(0,y;\mu) \end{bmatrix},
\label{eq:gslide}
\end{equation}
where $g_{\rm slide}$ is the sliding vector field.
The first component of \eqref{eq:gslide} is zero by the requirement of tangency to $x=0$.
This determines the value of $\theta$,
specifically $\theta = \frac{f_L}{f_L - f_R}$.
Upon substituting $\theta = \frac{f_L}{f_L - f_R}$ into the second component of \eqref{eq:gslide} we obtain
\begin{equation}
g_{\rm slide} = \frac{f_L g_R - f_R g_L}{f_L - f_R} \bigg|_{x=0}.
\label{eq:gslide2}
\end{equation}
In summary, on sliding regions orbits are governed by $\dot{y} = g_{\rm slide}(y;\mu)$,
where $g_{\rm slide}$ is given by \eqref{eq:gslide2}.

\subsection{Pseudo-equilibria}

Equilibria of $\dot{y} = g_{\rm slide}(y;\mu)$ are {\em pseudo-equilibria} of \eqref{eq:ode} and given by the roots of
\begin{equation}
h(y;\mu) = f_L(0,y;\mu) g_R(0,y;\mu) - f_R(0,y;\mu) g_L(0,y;\mu).
\label{eq:gslideNumerator}
\end{equation}
A pseudo-equilibrium is only exhibited by \eqref{eq:ode},
and said to be {\em admissible}, if it belongs to a sliding region.
Since \eqref{eq:ode} is piecewise-linear, $h(y;\mu)$ is a quadratic function of $y$.
Thus, generically, \eqref{eq:ode} has either no pseudo-equilibria, or two pseudo-equilibria.
Here we give conditions for the existence of admissible pseudo-equilibria.
This is proved below by directly calculating $h(y;\mu)$.

\begin{proposition}
Consider \eqref{eq:ode} with \eqref{eq:FJ} satisfying the assumptions of Theorem \ref{th:main}.
Let
\begin{align*}
c &= (a_{2L} b_{2R} - a_{2R} b_{2L}) \gamma, \\
d_L &= a_{2R}^2 \beta_L \,, \\
d_R &= a_{2L}^2 \beta_R \,,
\end{align*}
and $Q = c^2 - 2(d_L+d_R) c + (d_L-d_R)^2$.
\begin{enumerate}
\item
If $c \le |d_L - d_R|$ or $Q < 0$, then \eqref{eq:ode} has no admissible pseudo-equilibria for any $\mu \in \mathbb{R}$.
\item
If $c > |d_L - d_R|$ and $Q = 0$ [$Q > 0$], then \eqref{eq:ode} has one [two] admissible pseudo-equilibria for all $\mu \ne 0$.
\end{enumerate}
\label{pr:pseq}
\end{proposition}

\begin{figure}[t!]
\begin{center}
\setlength{\unitlength}{1cm}
\begin{picture}(9.5,17.2)
\put(0,11.8){\includegraphics[width=7.2cm]{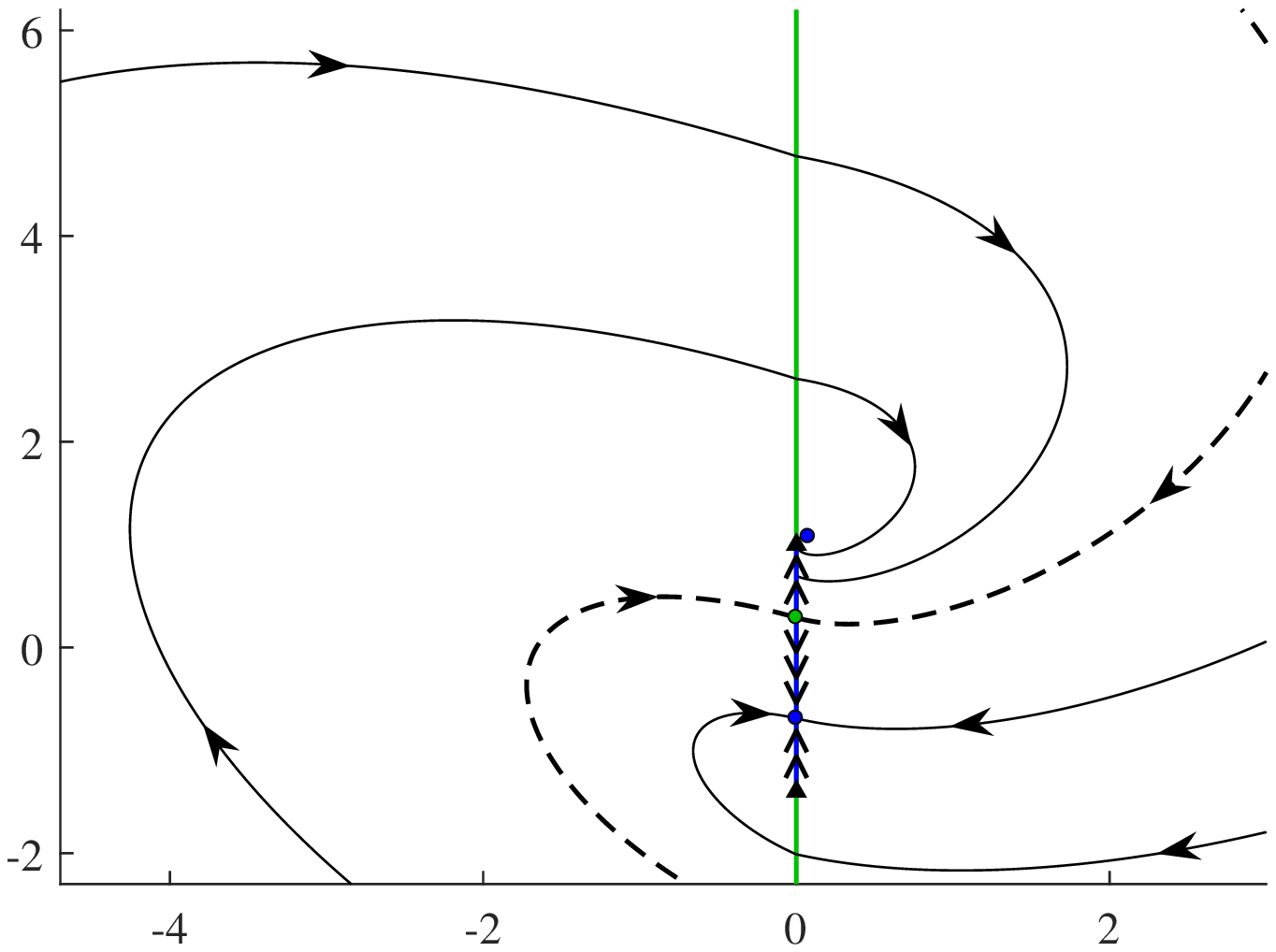}}
\put(0,5.9){\includegraphics[width=7.2cm]{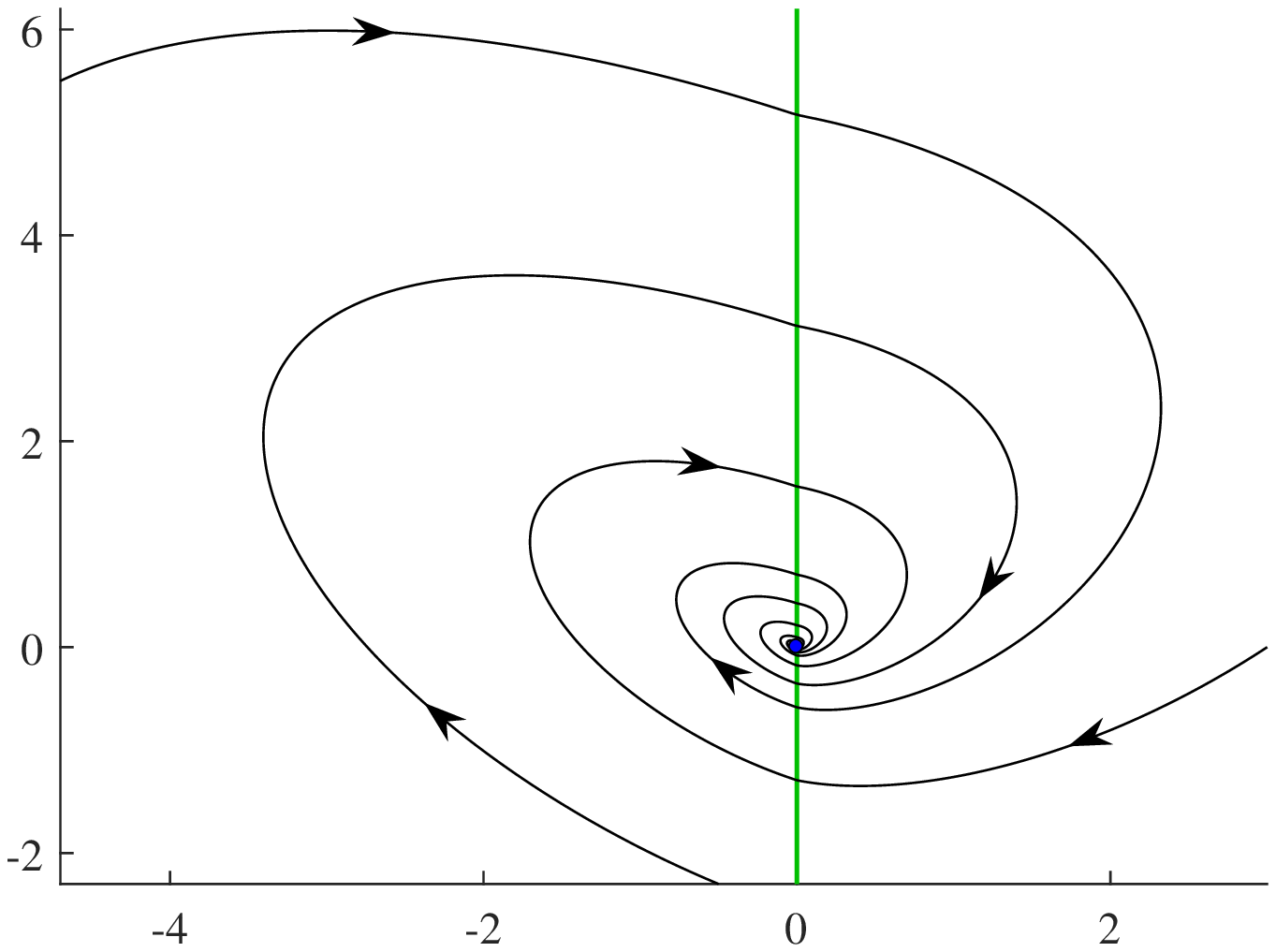}}
\put(0,0){\includegraphics[width=7.2cm]{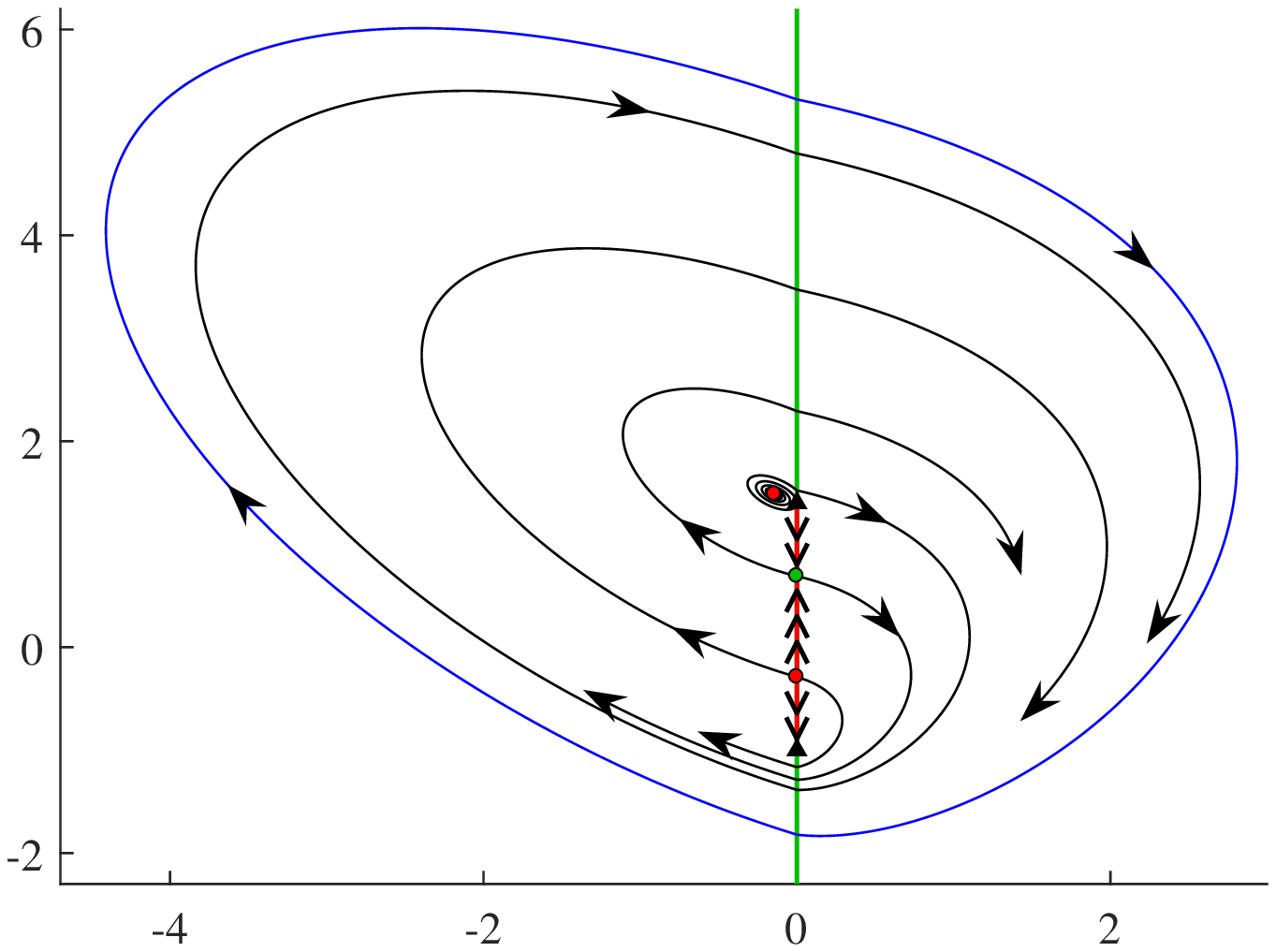}}
\put(4.56,11.8){\small $x$}
\put(0,14.81){\small $y$}
\put(7.7,14.81){\small {\bf a)}\hspace{2mm}$\mu = -1$}
\put(4.56,5.9){\small $x$}
\put(0,8.91){\small $y$}
\put(7.7,8.91){\small {\bf b)}\hspace{2mm}$\mu = 0$}
\put(4.56,0){\small $x$}
\put(0,3.01){\small $y$}
\put(7.7,3.01){\small {\bf c)}\hspace{2mm}$\mu = 1$}
\end{picture}
\caption{
Phase portraits of \eqref{eq:ex3}.
There exists a unique stable limit cycle for $\mu > 0$ by Theorem \ref{th:main},
and two pseudo-equilibria for all $\mu \ne 0$ by Proposition \ref{pr:pseq}.
In panel (a) the stable manifold of the saddle pseudo-equilibrium is dashed.
\label{fig:ppHLB4_3}
} 
\end{center}
\end{figure}

For the earlier example \eqref{eq:ex1}, we have $c = -2 \lambda_L < 0$.
Thus by Proposition \ref{pr:pseq}, \eqref{eq:ex1} has no
admissible pseudo-equilibria for all $\mu \in \mathbb{R}$.
Numerical investigations suggest that, for systems satisfying the assumptions of Theorem \ref{th:main},
the inequalities $c > |d_L - d_R|$ and $Q \ge 0$
are only satisfied in a relatively small fraction of parameter space.
A system that does satisfy these inequalities is
\begin{equation}
\begin{bmatrix} \dot{x} \\ \dot{y} \end{bmatrix} =
\begin{cases}
\begin{bmatrix}
\frac{3}{5} \,x + y - \frac{7}{5} \,\mu \\
-x - \frac{1}{2} \,y + \frac{3}{5} \,\mu
\end{bmatrix}, & x < 0, \\
\begin{bmatrix}
-x + y + \mu \\
-x - \frac{3}{10} \,y - \frac{2}{5} \,\mu
\end{bmatrix}, & x > 0.
\end{cases}
\label{eq:ex3}
\end{equation}
Here $a_{2L} = a_{2R} = 1$, $\beta_L = \beta_R = \frac{1}{10}$, and $\gamma = \frac{12}{5}$,
thus the conditions of Theorem \ref{th:main} are satisfied.
Also $c = \frac{12}{25}$ and $d_L = d_R = 0.1$,
thus $c > |d_L - d_R|$ and $Q = 0.0384 > 0$.
By Proposition \ref{pr:pseq}, \eqref{eq:ex3} has two admissible pseudo-equilibria for all $\mu \ne 0$.
These are shown in Fig.~\ref{fig:ppHLB4_3}.
When $\mu = -1$, one pseudo-equilibrium is stable.
The other is a saddle
and its stable manifold (dashed) forms the boundary between the basins of attraction of the stable pseudo-equilibrium
and the stable focus.
When $\mu = 1$, both pseudo-equilibria are unstable.
Here a stable limit cycle exists, by Theorem \ref{th:main}, because $\alpha \approx -0.63 < 0$.

\subsection{Proof of Proposition \ref{pr:pseq}}

Admissible pseudo-equilibria are roots of $h(y;\mu)$
that belong to the interval bounded by $\zeta_L$ and $\zeta_R$.
We first change variables so that this interval becomes $z \in [-1,1]$.
Let
\begin{equation}
\tilde{h}(z) = h \left( \frac{\zeta_L(\mu)+\zeta_R(\mu)}{2}
- \frac{\zeta_L(\mu)-\zeta_R(\mu)}{2} \,z; \mu \right),
\label{eq:tildeh}
\end{equation}
where by writing $\tilde{h}(z)$ we have suppressed the $\mu$-dependency for brevity.
Admissible pseudo-equilibria correspond to roots of $\tilde{h}(z)$ in $[-1,1]$.

\begin{figure}[t!]
\begin{center}
\setlength{\unitlength}{1cm}
\begin{picture}(7.2,5.4)
\put(0,0){\includegraphics[width=7.2cm]{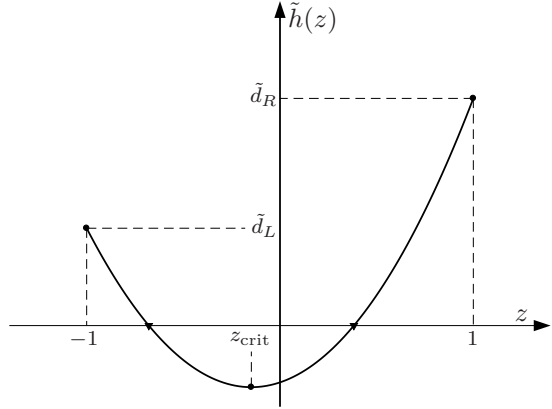}}
\put(6.72,1.14){\footnotesize $z$}
\put(3.71,5.05){\footnotesize $\tilde{h}(z)$}
\put(3.22,2.31){\scriptsize $\tilde{d}_L$}
\put(3.2,4.05){\scriptsize $\tilde{d}_R$}
\put(.8,.82){\scriptsize $-1$}
\put(2.93,.87){\scriptsize $z_{\rm crit}$}
\put(6.09,.82){\scriptsize $1$}
\end{picture}
\caption{
A sketch of the quadratic function $\tilde{h}(z)$, defined by \eqref{eq:tildeh}.
Roots of $\tilde{h}(z)$ in $[-1,1]$ correspond to admissible pseudo-equilibria of \eqref{eq:ode}.
\label{fig:quadraticSlidingFlow}
} 
\end{center}
\end{figure}

Since $\tilde{h}(z)$ is quadratic, it is completely determined by its second derivative
\begin{equation}
\tilde{c}
= \frac{d^2 \tilde{h}}{d z^2}
= \frac{(a_{2L} b_{2R} - a_{2R} b_{2L}) \gamma^2 \mu^2}{2 a_{2L}^2 a_{2R}^2}
= \frac{c \gamma \mu^2}{2 a_{2L}^2 a_{2R}^2},
\label{eq:tildec}
\end{equation}
and its values at $z = \pm 1$:
\begin{align}
\tilde{d}_L
&= \tilde{h}(-1)
= h \left( \zeta_L(\mu);\mu \right)
= \frac{\beta_L \gamma \mu^2}{a_{2L}^2}
= \frac{d_L \gamma \mu^2}{a_{2L}^2 a_{2R}^2},
\label{eq:tildedL} \\
\tilde{d}_R
&= \tilde{h}(1)
= h \left( \zeta_R(\mu);\mu \right)
= \frac{\beta_R \gamma \mu^2}{a_{2R}^2}
= \frac{d_R \gamma \mu^2}{a_{2L}^2 a_{2R}^2}.
\label{eq:tildedR}
\end{align}
That is, we can write $\tilde{h}(z)$ in terms of $\tilde{c}$, $\tilde{d}_L$, and $\tilde{d}_R$:
\begin{equation}
\tilde{h}(z) = \frac{1}{2} \left( \tilde{c} z^2 - \left( \tilde{d}_L - \tilde{d}_R \right) z
+ \tilde{d}_L + \tilde{d}_R - \tilde{c} \right).
\label{eq:tildeh2}
\end{equation}
Assuming $\tilde{c} \ne 0$, $\tilde{h}(z)$ has a unique critical point at
\begin{equation}
z_{\rm crit} = \frac{\tilde{d}_L - \tilde{d}_R}{2 \tilde{c}}.
\nonumber
\end{equation}
For the remainder of the proof we assume $\gamma > 0$
(the case $\gamma < 0$ can be dealt with similarly).
Choose any $\mu \ne 0$.
Then $\tilde{d}_L, \tilde{d}_R > 0$ (by \eqref{eq:tildedL}--\eqref{eq:tildedR},
because $\beta_L, \beta_R > 0$ by assumption).
Thus $\tilde{h}(z)$ has a root in $[-1,1]$ if and only if
(i) $\frac{d \tilde{h}}{d z}(-1) < 0$ and $\frac{d \tilde{h}}{d z}(1) > 0$, 
and (ii) $\tilde{h}(z_{\rm crit}) \le 0$, see Fig.~\ref{fig:quadraticSlidingFlow}.
Moreover, there is one root if $\tilde{h}(z_{\rm crit}) = 0$,
and two roots if $\tilde{h}(z_{\rm crit}) < 0$.

Condition (i) is $-\tilde{c} - \frac{\tilde{d}_L - \tilde{d}_R}{2} < 0$
and $\tilde{c} - \frac{\tilde{d}_L - \tilde{d}_R}{2} > 0$.
Since $\gamma > 0$, by \eqref{eq:tildec}--\eqref{eq:tildedR} these are equivalent to $c > |d_L - d_R|$.
Condition (ii) is $-\frac{(\tilde{d}_L - \tilde{d}_R)^2}{8 \tilde{c}}
+ \frac{\tilde{d}_L + \tilde{d}_R - \tilde{c}}{2} \le 0$.
If condition (i) holds, then $\tilde{c} > 0$ and so condition (ii) becomes
$\tilde{c}^2 - \left( \tilde{d}_L + \tilde{d}_R \right) \tilde{c}
+ \frac{1}{4} \left( \tilde{d}_L - \tilde{d}_R \right)^2 \ge 0$.
Since $\gamma > 0$, by \eqref{eq:tildec}--\eqref{eq:tildedR} this is equivalent to
$Q = c^2 - 2(d_L + d_R) c + (d_L - d_R)^2 \ge 0$.
\manualEndProof

\section{Three nested limit cycles}
\label{sec:threeLimitCycles}
\setcounter{equation}{0}

Here we provide an example to show that if the condition $a_{2L} \gamma \ge 0$
in Theorem \ref{th:main} is not satisfied, then the limit cycle may be non-unique.

Consider the system
\begin{equation}
\begin{bmatrix} \dot{x} \\ \dot{y} \end{bmatrix} =
\begin{cases}
\begin{bmatrix}
-\frac{4}{3} (x+\mu) + \frac{20}{3} \,y \\
-\frac{377}{750} (x+\mu) + \frac{26}{15} \,y
\end{bmatrix}, & x < 0, \\
\begin{bmatrix}
-\frac{19}{50} (x+\mu) + y \\
-(x+\mu) - \frac{19}{50} \,y
\end{bmatrix}, & x > 0,
\end{cases}
\label{eq:ex2}
\end{equation}
which is a simple transformation of the example given in \cite{BrMe13,HuYa12b}.
Here $a_{2L} = \frac{20}{3}$, $a_{2R} = 1$, $\beta_L = \frac{26}{25}$,
$\beta_R = \frac{2861}{2500}$, and $\gamma = -\frac{6}{5}$,
thus the conditions of Theorem \ref{th:main} are satisfied except $a_{2L} \gamma < 0$.

\begin{figure}[b!]
\begin{center}
\setlength{\unitlength}{1cm}
\begin{picture}(7.2,5.4)
\put(0,0){\includegraphics[width=7.2cm]{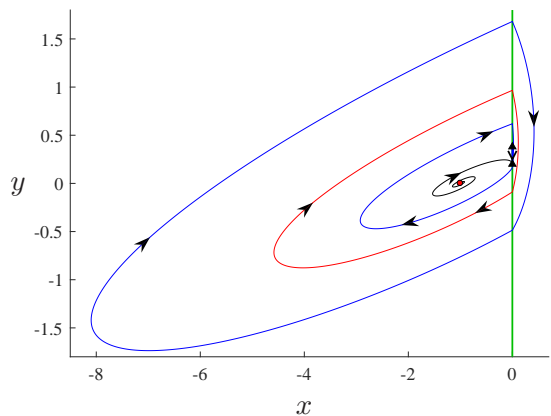}}
\put(3.8,0){\small $x$}
\put(0,2.98){\small $y$}
\end{picture}
\caption{
A phase portrait of \eqref{eq:ex2} with $\mu = 1$.
There exists an unstable focus, an attracting sliding region, and three nested limit cycles.
The middle limit cycle is unstable; the other two are stable.
\label{fig:ppHLB4_2}
} 
\end{center}
\end{figure}

\begin{figure}[b!]
\begin{center}
\setlength{\unitlength}{1cm}
\begin{picture}(7.2,3.85)
\put(0,0){\includegraphics[width=7.2cm]{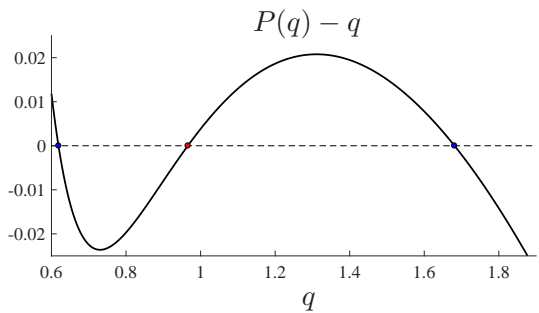}}
\put(4.05,0){\small $q$}
\put(3.4,3.65){\small $P(q) - q$}
\end{picture}
\caption{
The displacement $P(q)-q$, where $P$ is the Poincar\'e map of \eqref{eq:ex2} with $\mu = 1$
defined in the text.
Each fixed point $q^*$ of $P(q)$ corresponds to a limit cycle passing through $(x,y) = (0,q^*)$,
see Fig.~\ref{fig:ppHLB4_2}.
\label{fig:poinMapHLB4}
} 
\end{center}
\end{figure}

Since $a_{2L} \gamma < 0$, when $\mu = 1$ an attracting sliding region coexists
with the unstable focus, see Fig.~\ref{fig:ppHLB4_2}.
There are three limit cycles
and it is instructive to realise these as the fixed points of a Poincar\'e map.
Given $q > \zeta_R$, consider the forward orbit of $(x,y) = (0,q)$
and let $P(q)$ denote the $y$-value of the next intersection of this orbit with $x=0$ at a point with $y > \zeta_R$.
As seen in Fig.~\ref{fig:poinMapHLB4},
for sufficiently small values of $q$ we have $P(q) > q$ due to the unstable focus in $x < 0$,
whereas for sufficiently large values of $q$ we have $P(q) < q$ because $\alpha = -\frac{9}{50} < 0$.
Thus we expect $P(q)$ to have an odd number of fixed points --- indeed it has three
(established formally in \cite{BrMe13}).

It follows that \eqref{eq:ex2} has three limit cycles for all $\mu > 0$,
each with an amplitude proportional to $\mu$.
For $\mu < 0$, \eqref{eq:ex2} has no limit cycles.
Thus as the value of $\mu$ is increased from negative to positive,
three limit cycles are created as a stable focus
effectively turns into an unstable focus by colliding with $x=0$.

\section{Discussion}
\label{sec:conc}
\setcounter{equation}{0}

In this paper we have studied BEBs involving stable and unstable foci in Filippov systems.
The main result is a sufficient condition for a unique limit cycle to be created in the BEB, see Theorem \ref{th:main}.
The bifurcation resembles a Hopf bifurcation, and is one of $20$ Hopf-like bifurcations
of piecewise-smooth systems listed in \cite{Si18c}.
If the condition is not satisfied, three limit cycles may be created, see \S\ref{sec:threeLimitCycles}.
Up to two pseudo-equilibria can also arise, see \S\ref{sec:pseq}.

For simplicity only piecewise-linear systems have been considered
because, except in special cases, for general piecewise-smooth systems
the BEBs are expected to exhibit the same qualitative behaviour in a neighbourhood of the bifurcation.
As in the continuous setting \cite{SiMe07},
with nonlinear terms added to $F_L$ and $F_R$
hyperbolic equilibria, pseudo-equilibria, and limit cycles should be of a distance from the origin
that is asymptotically proportional to $|\mu|$, instead of directly proportional to $|\mu|$.
This will be treated carefully in \cite{Si18b}.

A major avenue for future research is to develop a comprehensive theory for
BEBs in systems of more than two dimensions.
Equilibria and pseudo-equilibria are well understood \cite{Si18d},
so it remains to characterise other invariant sets such as limit cycles and chaotic attractors.
As in the discrete-time case \cite{GlJe15},
a complete classification of BEBs in an arbitrary number of dimensions is surely infeasible,
but, as achieved in this paper, one can search for sufficient conditions for a BEB to behave in a certain fashion,
or determine when some form of dimension reduction is possible \cite{KoGl11,KuHo13,PrTe16}.

\end{document}